\documentclass{article}

\usepackage{amsfonts,amsmath,graphicx,subfigure}

\usepackage{bm}
\usepackage{geometry}
\usepackage{titling}
\usepackage{wrapfig}
\usepackage[backend=biber,style=authoryear-comp]{biblatex}
\usepackage{color}
\definecolor{darkgreen}{rgb}{0,0.8,0}

 \usepackage{authblk}
 \usepackage{titlesec}

\renewenvironment{abstract}
 {\small
 \par\noindent\textbf{\large\abstractname}\par\nobreak\smallskip}

\titleformat*{\section}{\large\bfseries}
\titleformat*{\subsection}{\normalsize\bfseries}
\titleformat*{\subsubsection}{\normalsize\bfseries}
\renewcommand\abstractname{ABSTRACT}

\newcommand\mycite[1]{\AtNextCite{\defcounter{maxnames}{1}}\citeauthor{#1}, \citeyear{#1}}
\newcommand{\citeauthorandyear}[2][]{(\mycite{#2})}

\title{A Novel Partitioned Approach for Reduced Order Model - Finite Element Model (ROM-FEM) and ROM-ROM Coupling}

\author[1]{Amy de Castro}
\author[2]{Paul Kuberry}
\author[3]{Irina Tezaur}
\author[4]{Pavel Bochev}
\affil[1]{Clemson University, agmurda@clemson.edu}
\affil[2]{Sandia National Laboratories, pakuber@sandia.gov}
\affil[3]{Sandia National Laboratories, ikalash@sandia.gov}
\affil[4]{Sandia National Laboratories, pbboche@sandia.gov}

\date{}

\begin{document}

\renewcommand{\refname}{REFERENCES}
\maketitle

\vspace{-0.5in}

\begin{abstract}

	Partitioned methods allow one to build a simulation capability for coupled problems 
	by reusing existing single-component codes.  
         In so doing, partitioned methods can shorten code development and validation times for multiphysics and multiscale applications. 
	In this work, we consider a scenario in which one or more of the ``codes'' being coupled are projection-based reduced order models (ROMs),
	introduced to lower the computational cost associated with a particular component.
%
	We simulate this scenario by considering a model interface problem that is discretized independently on two non-overlapping subdomains. We then formulate a partitioned scheme for this problem that allows the coupling between a ROM ``code'' for one of the subdomains with a finite element 
	model (FEM) or ROM ``code'' for the other subdomain.
	The  ROM ``codes'' are constructed by performing proper orthogonal decomposition (POD) 
	on a snapshot ensemble
	to obtain a low-dimensional reduced order basis, followed by a Galerkin projection onto this basis.
The ROM and/or FEM ``codes'' on each subdomain are then coupled using a Lagrange multiplier
representing the interface flux. To partition the resulting monolithic problem, we first eliminate the flux through a dual Schur complement.
Application of an explicit time integration scheme to the transformed monolithic problem
decouples the subdomain equations, allowing their independent
solution for the next time step. We show numerical results that demonstrate the proposed method’s
efficacy in achieving both ROM-FEM and ROM-ROM coupling.

\end{abstract}

\section{INTRODUCTION} \label{AdC:sec:intro}
Partitioned schemes enable the rapid development of simulation capabilities for coupled problems from 
existing codes for the individual sub-models; see, e.g., \citeauthorandyear{deBoer_07_CMAME} for examples. Besides being a cost-effective alternative to the development of monolithic multiphysics codes from scratch, a partitioned approach can also improve simulation efficiency by employing codes tailored to the salient physics characteristics of the sub-models.

Typically, the sub-model codes in partitioned schemes implement high-fidelity full-order  models (FOMs) based on conventional discretizations such as finite elements, finite volumes or finite differences. However, it is not uncommon to encounter situations in which one or more of these full order models become performance bottlenecks. For example, in blast-on-structure simulations \citeauthorandyear{Bessette_03a_INPROC}, calculation of the wave propagation by a high-fidelity scheme can be computationally expensive and is often replaced with direct structure loading by means of simplified boundary conditions  derived using analytic techniques \citeauthorandyear{Randers-Pehrson_97_ARL}. However, such conditions assume simple geometries and cannot account for wave interactions with more complex fluid-structure interfaces. A better alternative in this context would be a \emph{hybrid} partitioned scheme in which the expensive full-order sub-model is replaced by a computationally efficient, yet physically faithful, \emph{reduced order model} (ROM).

To demonstrate the potential of a coupling approach of the type described above, we formulate herein a new hybrid explicit partitioned scheme that enables the coupling of conventional finite element models (FEM) with projection-based ROMs (more specifically, ROMs constructed using the Proper Orthogonal Decomposition (POD)/Galerkin projection method (\mycite{Holmes:1988}; \mycite{Holmes:1996}; \mycite{Sirovich:1987}). We describe and develop our methodology in the context of a generic
advection-diffusion transmission problem posed on a  decomposition of the physical domain into two non-overlapping subdomains. Although simple 
and comprised of a single physics, this problem configuration is sufficient to simulate a typical setting for the development of a partitioned scheme.

Our scheme extends the approach in \citeauthorandyear{AdC:CAMWA}, which starts from a monolithic 
formulation of the transmission problem, uses a Schur complement to obtain an approximation of 
the interface flux, and then inserts this flux as a Neumann boundary condition into each subdomain problem. 
Application of an explicit time integration scheme to this transformed monolithic problem decouples its subdomain problems and allows their independent solution.


In addition to enabling a hybrid partitioned analysis for coupled problems, our approach can also be used to perform a \emph{hybrid} reduced order model - full order model (ROM-FOM) analysis (\mycite{Lucia:2001}; \mycite{Lucia:2003}; \mycite{LeGresley:2003};
\mycite{LeGresley:2005}; \mycite{Buffoni:2007}; \mycite{Baiges:2013}; \mycite{Corigliano:2015}). In this approach, the physical domain of a given, usually single physics, partial differential equation (PDE) problem is decomposed into two or more subdomains, and either a ROM or a FOM is constructed in each subdomain based on the solution characteristics. The resulting models are then coupled in some way to obtain a global solution on the physical domain in its entirety. Such an analysis can mitigate robustness and accuracy issues of projection-based model order reduction, especially when applied to highly non-linear and/or convection-dominated problems. 

In contrast to traditional partitioned schemes (\mycite{Gatzhammer_14_THESIS};
\mycite{Piperno_01_CMAME}; \mycite{Banks_17_JCP}) and methods for hybrid ROM-FOM analyses (\mycite{LeGresley:2003}; \mycite{LeGresley:2005};
\mycite{Buffoni:2007}; \mycite{Cinquegrana:2011};  \mycite{Maier:2014}),  our framework is monolithic rather than
iterative, enabling one to obtain the coupled ROM-ROM or ROM-FEM solution in a single shot.  
Also, unlike the work in (\mycite{Ammar:2011}; \mycite{Iapichino:2016};
 \mycite{Eftang:2013}; \mycite{Eftang:2014}; \mycite{Hoang:2021}), 
there is no need in our formulation to construct boundary, port, or skeleton bases for enforcing
inter-subdomain compatibility.  Furthermore, while our formulation shares some commonalities with existing Lagrange multiplier-based coupling methods such as those of  (\mycite{Lucia:2001}; \mycite{Lucia:2003};
 \mycite{Maday:2004}; \mycite{Antil:2010};
\mycite{Corigliano:2013}; \mycite{Corigliano:2015};
 \mycite{Kerfriden:2013};  \mycite{Radermacher:2014}; 
 \mycite{Baiges:2013}), we emphasize that our approach is fundamentally different
from these methods in that it enables the complete decoupling of the underlying models (ROMs and/or FOMs) at each
time-step of the time-integration scheme used to advance the discretized PDE forward in time.  Importantly, our 
methodology delivers a smooth and accurate solution without the need to introduce \textit{ad hoc} correction/stabilization
terms, such as those proposed in (\mycite{LeGresley:2003}; \mycite{LeGresley:2005}; \mycite{Baiges:2013}).

The remainder of this paper is organized as follows. In Section \ref{sec:interface}, we introduce our model transmission problem, derive the relevant monolithic formulation and discretize it in space.  
Section  \ref{AdC:FEM scheme} explains the elimination of the Lagrange multiplier through a dual Schur complement, which transforms the semi-discrete system into another coupled problem that serves as a basis for the development of our partitioned scheme.
%
In Section \ref{AdC:sec:partitioned}, we describe our POD/Galerkin ROM construction
methodology, and detail the application of the approach described in Section \ref{AdC:FEM scheme} to ROM-FEM and ROM-ROM coupling.  
We evaluate the performance of the proposed scheme on a two-dimensional (2D) model problem in Section \ref{sec:numerical}.  Finally, conclusions are offered in  Section \ref{sec:conc}.

\section{A MODEL TRANSMISSION PROBLEM} \label{sec:interface}

\begin{figure}[!htbp]
\begin{center}
	{\includegraphics[scale=.55]{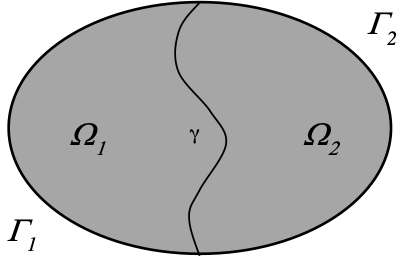}} \caption{Non-overlapping subdomains} \label{fig:genDomain}
\end{center}\end{figure}

We consider a bounded region $\Omega\subset \mathbb{R}^d$,  $d = 2,3$ with a Lipschitz-continuous boundary $\Gamma$. We assume that $\Omega$ is divided into two non-overlapping subdomains $\Omega_1$ and  $\Omega_2$, each with boundary 
$\partial \Omega_i$ for $i=1,2$.  Let $\gamma$ denote the interface shared between the two subdomains, and let $\Gamma_i = \partial\Omega_i \backslash \gamma$ for $i = 1,2$, as illustrated in Figure \ref{fig:genDomain}.
We take $\mathbf{n}_\gamma$ to be the unit normal on the interface pointing toward $\Omega_2$. 
We use a setting comprising two non-overlapping domains to avoid technical complications that are not germane to the core topic of the paper. Our approach can be extended to configurations involving multiple domains as long as one incorporates a proper mechanism to handle floating subdomains, such as the techniques in \citeauthorandyear{Bochev_05_SIREV}.

We consider a model transmission problem given by the advection-diffusion equation:

\begin{align}\label{AdC:strongForm}
\begin{split}
\dot{\varphi_i} - \nabla \cdot F_i (\varphi_i) &= f_i \hspace{5mm} \text{ on } \Omega_i \times [0,T] \\
\varphi_i &= g_i \hspace{5mm} \text{ on } \Gamma_i \times [0,T], \hspace{2mm} i=1,2, 
\end{split}
\end{align}
where the over-dot notation denotes differentiation in time, 
	the unknown $\varphi_i$ is a scalar field, $F_i(\varphi_i) = \kappa_i \nabla \varphi_i - \mathbf{u} \varphi_i$ is the total flux function, $\kappa_i > 0$ is the diffusion coefficient in $\Omega_i$, and $\mathbf{u}$ the velocity field.
We augment \eqref{AdC:strongForm} with initial conditions:
\begin{align}\label{AdC:masterics}
	\varphi_i(\mathbf{x},0) = \varphi_{i,0}(\mathbf{x}) \hspace{5mm} \text{ in } \Omega_i, \hspace{2mm} i=1,2\,.
\end{align}
Along the interface $\gamma$, we enforce continuity of the states and continuity of the total flux, giving rise 
	to the following interface conditions:
\begin{align}\label{AdC:interfaceMethod1}
\varphi_1(\mathbf{x},t) - \varphi_2(\mathbf{x},t) = 0 \hspace{2mm} \text{ and } F_1(\mathbf{x},t) \cdot \mathbf{n}_\gamma = F_2(\mathbf{x},t) \cdot \mathbf{n}_\gamma \hspace{2mm} \text { on } \gamma \times [0,T].
\end{align}
We note that one also has the option to enforce only  equilibrium of  the diffusive flux exchanged between the two subdomains. We do not consider this option here, as the resulting partitioned scheme will be similar to the one obtained by enforcing continuity of the total flux.

	In contrast to conventional, loosely coupled partitioned schemes (see, e.g., \citeauthorandyear{deBoer_07_CMAME}), our approach starts from a well-posed monolithic formulation of \eqref{AdC:strongForm}--\eqref{AdC:interfaceMethod1}. 
To obtain this formulation let $V:=H^1_\Gamma(\Omega_1) \times H^1_\Gamma(\Omega_2) \times H^{-1/2}(\gamma)$. Using a Lagrange multiplier to enforce continuity of states, i.e., the first condition in \eqref{AdC:interfaceMethod1}, yields the following monolithic weak  problem:
\emph{find} $\{\varphi_1, \varphi_2, \lambda\} \in C^1([0,T];V)$, \emph{such that for all $t\in (0,T]$}
\begin{align}
\begin{split}\label{AdC:overallM1}
(\dot{\varphi_1}, \nu)_{\Omega_1} + (\kappa_1 \nabla \varphi_1, \nabla \nu)_{\Omega_1} - (\mathbf{u} \varphi_1, \nabla \nu)_{\Omega_1} + (\lambda, \nu)_\gamma &= (f_1, \nu)_{\Omega_1}    \hspace{5mm} \forall \nu \in H^1_\Gamma(\Omega_1) \\
(\dot{\varphi_2}, \eta)_{\Omega_2} + (\kappa_2 \nabla \varphi_2, \nabla \eta)_{\Omega_2} - (\mathbf{u} \varphi_2, \nabla \eta)_{\Omega_2} - (\lambda, \eta)_\gamma  &= (f_2, \eta)_{\Omega_2}    \hspace{5mm} \forall \eta \in H^1_\Gamma(\Omega_2) \\
(\varphi_1, \mu)_\gamma - (\varphi_2, \mu)_\gamma &= 0 \hspace{17mm} \forall \mu \in H^{-1/2}(\gamma).
\end{split}
\end{align}
It is easy to see that the Lagrange multiplier $\lambda$ is the flux exchanged through the interface, i.e., $\lambda=F_1 \cdot \mathbf{n}_\gamma = F_2 \cdot \mathbf{n}_\gamma$. This observation is at the core of our partitioned method formulation. Indeed, if we could somehow determine $\lambda$, then each subdomain problem becomes a well-posed mixed boundary value problem with a Neumann condition on $\gamma$ provided by $\lambda$:
\begin{align}
\begin{split}\label{AdC:govM1}
\dot{\varphi_i} - \nabla \cdot F_i(\varphi_i) &= f_i \hspace{15mm} \text{ on } \Omega_i \times [0,T] \\
\varphi_i &= g_i \hspace{15mm} \text{ on } \Gamma_i \times [0,T] \\
F_i(\varphi_i) \cdot \mathbf{n}_i &= (-1)^i \lambda \hspace{8mm} \text{ on } \gamma \times [0,T]
\end{split}
\,,\quad i=1,2 \, .
\end{align}
In other words, knowing $\lambda$ could allow us to decouple the subdomain equations and solve them independently.
Of course, this cannot be done within the framework of \eqref{AdC:overallM1}, which is a fully coupled problem in terms of the states $\phi_i$ \emph{and} the Lagrange multiplier $\lambda$. However, an independent estimation of $\lambda$ may be possible in the context of a discretized version of this coupled problem.

\subsection{A SEMI-DISCRETE MONOLITHIC FORMULATION}\label{AdC:sec:semi}
Let $V^h\subset V$ be a conforming finite element space spanned by a basis $\{\nu_{i}, \eta_{j}, \mu_k\}$; $i=1,\ldots, N_1$; $j=1,\ldots, N_2$; $k=1,\ldots, N_{\gamma}$.  A finite element discretization of \eqref{AdC:overallM1} yields the following system of Differential Algebraic Equations (DAEs):
\begin{align}\label{AdC:matrixfull-index2}
\begin{split}
M_1 \dot{\bm{\Phi}}_1 + G_1^T \boldsymbol{\lambda} &= \mathbf{\overline{f}}_1(\bm{\Phi}_1)  \\
M_2 \dot{\bm{\Phi}}_2 - G_2^T \boldsymbol{\lambda} &= \mathbf{\overline{f}}_2(\bm{\Phi}_2) \\
	G_1 \bm{\Phi}_1 - G_2 \boldsymbol{{\Phi}}_2 &= \boldsymbol{0},
\end{split}
\end{align}
where for $r=1,2$, $\bm{\Phi}_r$ are the coefficient vectors corresponding to $\varphi_r$, $M_r$ are the mass matrices, the right hand side vector
$\mathbf{\overline{f}}_r(\bm{\Phi}_r) := \mathbf{f}_r - (D_r + A_r) \bm{\Phi}_r$ with $D_r, A_r$ corresponding to the diffusive and advective flux terms, respectively, and $G_r$ are the matrices enforcing the (weak) continuity of the states.
Assembly of these matrices is standard, for example, $(M_1)_{ij} = (\nu_j,\nu_i)_{\Omega_1}$, $(D_2)_{ij} = \kappa_2 (\nabla \eta_j, \nabla \eta_i)_{\Omega_2}$; $(G_1)_{i,j} = (\nu_j, \mu_i)_{\gamma}$;  $(G_2)_{i,j} = (\eta_j, \mu_i)_\gamma $, and so on. 
We note here that the space for the Lagrange multiplier $\lambda$ can be taken to be the trace of the finite element space on either of $\Omega_1$ or $\Omega_2$; either choice will be stable. In practice, using the coarser of the two interface spaces for the Lagrange multiplier space improves accuracy; see \citeauthorandyear{AdC:CAMWA} and \citeauthorandyear{AdC:RINAM} for details and discussion.

%

\section{EXPLICIT PARTITIONED SCHEME FOR FEM-FEM COUPLING}\label{AdC:FEM scheme}

In this section, we briefly review the Implicit Value Recovery (IVR) scheme \citeauthorandyear{AdC:CAMWA}, which provides the basis for our new hybrid partitioned approach. Then, in Section  \ref{AdC:sec:partitioned}, we discuss extensions of IVR to include a ROM in one or both subdomains.

 The IVR scheme \citeauthorandyear{AdC:CAMWA} is predicated on the ability to express $\bm{\lambda}$ as an implicit function of the subdomain states. This, however, is not possible for \eqref{AdC:matrixfull-index2} because it is an Index-2 Hessenberg DAE. 
In \citeauthorandyear{AdC:CAMWA}, we resolved this issue by differentiating the constraint equation in time. This step reduced the index of 
 \eqref{AdC:matrixfull-index2} and produced the following Index-1 Hessenberg DAE:
\begin{align}\label{AdC:matrixfull}
\begin{split}
M_1 \dot{\bm{\Phi}}_1 + G_1^T \boldsymbol{\lambda} &= \mathbf{\overline{f}}_1(\bm{\Phi}_1)  \\
M_2 \dot{\bm{\Phi}}_2 - G_2^T \boldsymbol{\lambda} &= \mathbf{\overline{f}}_2(\bm{\Phi}_2) \\
G_1 \dot{\bm{\Phi}}_1 - G_2 \dot{\boldsymbol{{\Phi}}}_2 &= 0\;\;.
\end{split}
\end{align}
Assuming  the initial data are continuous across $\gamma$, the new constraint $(\dot{\varphi_1}, \mu)_\gamma - (\dot{\varphi_2}, \mu)_\gamma = 0$ is equivalent to the original one, i.e., \eqref{AdC:matrixfull} is equivalent to the original monolithic problem \eqref{AdC:matrixfull-index2}. In what follows we refer to \eqref{AdC:matrixfull} as the FEM-FEM model. This model can be written in matrix form as:
\begin{equation}\label{AdC:coupled425}
\begin{bmatrix}
M_1 & 0 & G_1^T \\ 0 & M_2 & -G_2^T \\ G_1 & -G_2 & 0
\end{bmatrix}
\begin{bmatrix}
\dot{\bm{\Phi}}_1 \\ \dot{\bm{\Phi}}_2 \\ \boldsymbol{\lambda}
\end{bmatrix}
= \begin{bmatrix}
\mathbf{\overline{f}}_1(\bm{\Phi}_1) \\ \mathbf{\overline{f}}_2(\bm{\Phi}_2) \\ 0
\end{bmatrix}\,.
\end{equation}
%
%
To explain IVR it is further convenient to write \eqref{AdC:coupled425} in the canonical semi-explicit DAE form:
\begin{align}\label{AdC:algDAE}
\begin{split}
\dot{y} &= f(t,y,z) \\
0 &= g(t,y,z)
\end{split}
\end{align}
where $y = (\bm{\Phi}_1, \bm{\Phi}_2)$ is the differential variable, $z = \boldsymbol{\lambda}$ is the algebraic variable,
\begin{align}
\begin{split}
f(t,y,z) &= \begin{pmatrix}
M_1^{-1} \Big(\overline{\mathbf{f}}_1(\bm{\Phi}_1) - G_1^T \boldsymbol{\lambda} \Big) \\
M_2^{-1} \Big( (\overline{\mathbf{f}}_2(\bm{\Phi}_2) + G_2^T \boldsymbol{\lambda} \Big)
\end{pmatrix}
\end{split}\;,
\end{align}
and
\begin{align}\label{AdC:defG}
g(t,y,z) = S \boldsymbol{\lambda} - G_1 M_1^{-1} \overline{\mathbf{f}}_1(\bm{\Phi}_1) + G_2 M_2^{-1} \overline{\mathbf{f}}_2(\bm{\Phi}_2)\,.
\end{align}
The matrix $S := G_1 M_1^{-1} G_1^T + G_2 M_2^{-1} G_2^T$ in \eqref{AdC:defG} is the Schur complement of the upper left $2\times 2$ block submatrix of the matrix in \eqref{AdC:coupled425}.  
It can be shown that the Schur complement $S$ is nonsingular; see Proposition 4.1 in \citeauthorandyear{AdC:CAMWA}. This implies that the Jacobian $\partial_z g = S$ is also nonsingular for all $t$.
As a result, the equation $g(t,y,z) = 0$ defines $z$ as an implicit function of the differential variable. After solving this equation for the algebraic variable and inserting the solution $\bm{\lambda}(\bm{\Phi}_1,\bm{\Phi}_2)$  into \eqref{AdC:coupled425} we obtain a coupled system of ODEs in terms of the states:
\begin{equation}\label{AdC:decoupled}
\begin{bmatrix}
M_1 & 0 \\ 0 & M_2 
\end{bmatrix}
\begin{bmatrix}
\dot{\bm{\Phi}}_1 \\ \dot{\bm{\Phi}}_2 
\end{bmatrix}
= \begin{bmatrix}
\mathbf{\overline{f}}_1(\bm{\Phi}_1) - G_1^T \bm{\lambda}(\bm{\Phi}_1,\bm{\Phi}_2) \\ \mathbf{\overline{f}}_2(\bm{\Phi}_2) + G_2^T \bm{\lambda}(\bm{\Phi}_1,\bm{\Phi}_2) 
\end{bmatrix}\;.
\end{equation}
The IVR scheme is based on the observation that application of an explicit time integration scheme to discretize \eqref{AdC:decoupled} in time effectively decouples the equations and makes it possible to solve them independently.

The IVR algorithm for solving the coupled system is now as follows. Let $D^n_t(\bm{\Phi})$ be a forward time differencing operator such as the Forward Euler operator $D^n_t(\bm{\Phi}) = (\bm{\Phi}^{n+1} - \bm{\Phi}^n ) / \Delta t$. For each time step $t^n$:
\begin{enumerate}
\item \textit{Compute modified forces}: for $i=1,2$ use $\bm{\Phi}_i^n$ to compute the vector \begin{equation*} \widetilde{\boldsymbol{f}}_{i}^{n} :=\overline{\mathbf{f}}_i(\bm{\Phi}_i^n) = \mathbf{f}_i - (D_i + A_i)  \bm{\Phi}_i^n. \end{equation*}
\item \textit{Compute the Lagrange multiplier}: solve the Schur complement system 
\begin{equation*}
\Big( G_1 M_1^{-1} G_1^T + G_2 M_2^{-1} G_2^T \Big) \bm{\lambda}^n = G_1 M_1^{-1} \widetilde{\bm{f}}_1^n - G_2 M_2^{-1} \widetilde{\boldsymbol{f}_2^n}
\end{equation*}
 for $\bm{\lambda}^n$. Compute $G_1^T \bm{\lambda}^n$ and $G_2^T \bm{\lambda}^n$.
 \item \textit{Update the state variables:} for $i=1,2$, solve the systems 
 \begin{equation*} M_i D^n_t(\boldsymbol{\Phi}_i) = \widetilde{\boldsymbol{f}}_{i}^{n} + (-1^i) G_i^T \bm{\lambda}^n. \end{equation*}
\end{enumerate}

\section{DEVELOPMENT OF HYBRID PARTITIONED SCHEMES}\label{AdC:sec:partitioned}

In this section, we present the extension of the IVR method, described in Section \ref{AdC:FEM scheme}, to a \emph{hybrid partitioned scheme}, which couples ROM to FEM or to another ROM. Specifically, in Section \ref{AdC:sec:ROMsystem}, we present the details for the case where a projection-based ROM is employed in one of the two subdomains
shown in Figure \ref{fig:genDomain}; then, in Section \ref{sec:ROM-ROM}, we briefly describe the ROM-ROM extension of our coupling methodology.

To define the ROM component within our coupling, we use a proper orthogonal decomposition (POD) approach. A typical POD-based model order reduction 
is comprised of two distinct stages. In the first stage, one uses samples obtained by solving a suitable FOM to construct a reduced basis, usually by computing a truncated singular value decomposition (SVD) of the sample set. We discuss this stage in Section \ref{ROMsec}. At the second stage, one replaces the conventional finite element test and trial functions in a weak formulation of the governing equations by reduced basis functions.  This stage projects the weak problem onto the reduced basis and is discussed in Section \ref{AdC:sec:ROMsystem}.

Obtaining a quality ROM that is both computationally efficient and  accurate in the predictive regime is a non-trivial endeavor on its own. Since our main goal is the development of the hybrid partitioned approach rather than the ROM, in this work we will follow standard, established procedures to obtain the necessary ROMs.

\subsection{REDUCED BASIS CONSTRUCTION}\label{ROMsec}

Without loss of generality, we shall describe the first stage of the POD-based model order reduction for $\Omega_1$. In this work,
we have adopted a workflow in which the snapshots in $\Omega_1$ are collected by performing a global (uncoupled) FEM simulation in $\Omega$, restricting the resulting finite element solution to  $\Omega_1$, and then sampling the restricted solution over $m$ uniform time steps. Let $\Delta_s t$ denote the sampling time step, $t_k = k (\Delta_s t)$, $k=1, \ldots m$, the sampling time points, and $\bm{\Phi}_1(t_k)\in\mathbb{R}^{N_1}$ 
the $k$th snapshot, i.e., the coefficient vector of the restricted finite element solution at $t_k$.

We arrange the snapshots in an $N_1 \times m$ matrix $X$ whose $k$th column is the $k$th snapshot $\bm{\Phi}_1(t_k)$. The coefficients in each snapshot form two distinct groups. The first one contains the coefficients associated with the nodes on the Dirichlet boundary $\Gamma_1$.  These coefficients contain the values of the boundary condition function $g_1$ at the these nodes, and so we call them \emph{Dirichlet coefficients}. 
The coefficients in the second group correspond to the nodes in the interior of $\Omega_1$ and the nodes on the interface $\gamma$. 
We refer to these coefficients as the \emph{free coefficients}, as they are the unknowns in the finite element discretization of the subdomain PDE \eqref{AdC:govM1} on $\Omega_1$.

Performing a POD-based model order reduction for problems with Dirichlet boundary conditions requires some care in the generation of the reduced basis and the subsequent imposition of the Dirichlet conditions on the ROM solution. Herein, we use an approach that represents an extension of a common finite element technique that imposes essential boundary conditions via a boundary interpolant of the data $g_1$; see \citeauthorandyear{Gunzburger:2007} for more details. Below we describe how this technique is applied to the generation of the reduced basis, and in Section \ref{AdC:sec:ROMsystem}, we explain the imposition of the boundary conditions within the ROM formulation.

Let $\bm{\beta}_k \in \mathbb{R}^{N_1}$  denote a vector whose free coefficients are all set to zero and whose Dirichlet coefficients are set to the nodal values of the boundary data at $t_k$, that is, 
\begin{equation}\label{eq:beta}
(\bm{\beta}_k)_i  = 
\left\{\begin{array}{rl}
g_1(\mathbf{x}_i,t_k) & \mbox{if  $\mathbf{x}_i\in\Gamma_1$} \\[1ex]
0         & \mbox{if  $\mathbf{x}_i\in\Omega_1\cup\gamma$}  
\end{array}\;.
\right.
\end{equation}
Following  \citeauthorandyear{Gunzburger:2007} we define the \emph{adjusted} snapshot matrix  $X_0$ by subtracting\footnote{Note that the net effect of this computation zeros out the rows in $X$ containing  Dirichlet coefficients while leaving the rows containing free coefficients unchanged. Thus in practice, one may manually zero the Dirichlet rows of $X$. } $\bm{\beta}_k$ from the $k$th column of $X$, i.e., we set the $k$th column of $X_0$ to $\bm{\Phi}_1(t_k) - \bm{\beta}_k$.

Next, we compute the singular value decomposition of the adjusted matrix, $X_0 = U_0 \Sigma_0 V^T_0$, and choose an integer $N_R \ll N_1$. The reduced basis is then defined as the first $N_R$ left singular vectors of the SVD decomposition, i.e., the first $N_R$ columns of $U_0$. We denote the matrix  containing these columns by $\widetilde{U}_0$.

Each column of $\widetilde{U}_0$ can be mapped to a finite element function whose nodal coefficients are the entries in this column. These finite element functions can be construed as a new reduced order basis for the finite element space. Note that these basis functions  are globally supported rather than locally supported, as is the case with traditional finite element basis functions. Thus, using the reduced basis functions as test and trial functions in a weak formulation of \eqref{AdC:govM1}  results in dense algebraic problems.  Consequently, an effective ROM requires $N_R$ to be as small as possible. A simple approach is to choose a tolerance level $\delta$ and remove the columns of $U_0$ corresponding to all singular values that are less than $\delta$. We note that $\delta$ should be such that no columns of $U_0$ are retained which correspond to singular values sufficiently close to 0. These columns of $U_0$ span a near null space, which we do not want to retain as part of the reduced basis.

\subsection{IVR EXTENSION TO ROM-FEM COUPLING}\label{AdC:sec:ROMsystem}
To extend the IVR scheme in Section \ref{AdC:FEM scheme} from a FEM-FEM to a ROM-FEM coupling with a ROM on $\Omega_1$, we will perform the second model order reduction stage directly in the monolithic formulation of the model problem. 
Formally, this amounts to discretizing the first equation in \eqref{AdC:overallM1} using the 
global basis functions (POD modes)
corresponding to the columns of $\widetilde{U}_0$ instead of the standard finite element basis functions. In practice, for linear problems,
the matrices defining the ROM can be be easily obtained from the already assembled full order model matrices. Thus, we will implement the second stage using the transformed semi-discrete monolithic problem, i.e., the Index-1 DAE \eqref{AdC:matrixfull}.

For simplicity, in discussing this stage, we shall assume that the Dirichlet boundary condition function $g_1$ is independent of time. In this case,  the vectors $\bm{\beta}_k$ defined in \eqref{eq:beta} are identical to a vector $\bm{\beta}$ whose free coefficients are zero and Dirichlet coefficients are the nodal values of $g_1$. To obtain the ROM on $\Omega_1$ we perform a state transformation of the first equation in \eqref{AdC:matrixfull} by inserting the ansatz $\bm{\Phi}_1 = \widetilde{U}_0 \bm{\varphi}_R + \bm{\beta}$ into that equation. Then, we multiply the first equation by $\widetilde{U}^T_0$ to obtain the following ROM-FEM monolithic problem:
\begin{align}\label{AdC:ROMsystem}
\begin{split}
\widetilde{M}_1\dot{\boldsymbol{\varphi}}_R + \widetilde{G}^T_1 \boldsymbol{\lambda} 
&= \widetilde{U}^T_0 \overline{\mathbf{f}}_1(\widetilde{U}_0 \bm{\varphi}_R  + \bm{\beta})\\
M_2 \dot{\bm{\Phi}}_2 - G_2^T \boldsymbol{\lambda} &= \overline{\mathbf{f}}_2(\bm{\Phi}_2) \\
\widetilde{G}_1 \dot{\boldsymbol{\varphi}}_R - G_2 \dot{\bm{\Phi}}_2 &= 0,
%
\end{split}
\end{align}
where $\widetilde{M}_1: = \widetilde{U}^T_0 M_1  \widetilde{U}_0 $ and $\widetilde{G}_1^T := \widetilde{U}^T_0 G_1^T $.
Note that the first equation is now of size $N_R$. Let $y = (\bm{\varphi}_R, \bm{\Phi}_2)$ be the differential variable, and $z = \boldsymbol{\lambda}$ the algebraic variable. As in Section \ref{AdC:FEM scheme}, the ROM-FEM monolithic system \eqref{AdC:ROMsystem} is an index-1 DAE having the same canonical form as \eqref{AdC:algDAE} but with:
\begin{align}
\begin{split}
f(t,y,z) 
&= \begin{pmatrix}
\widetilde{M}_1^{-1} \Big( \widetilde{U}^T_0\overline{\mathbf{f}}_1(\widetilde{U}_0 \bm{\varphi}_R  + \bm{\beta}) - 
	\widetilde{G}_1^T \boldsymbol{\lambda} \Big) \\[2ex]
M_2^{-1} \Big( \overline{\mathbf{f}}_2(\bm{\Phi}_2) + G_2^T \boldsymbol{\lambda} \Big),
\end{pmatrix}
\end{split}
\end{align}
and
\begin{align}
g(t,y,z) = \widetilde{S} \boldsymbol{\lambda} 
-  \widetilde{G}_1\widetilde{M}_1^{-1}
\Big( \widetilde{U}^T_0 \overline{\mathbf{f}}_1(\widetilde{U}_0 \bm{\varphi}_R  + \bm{\beta})\Big) + G_2 M_2^{-1} \overline{\mathbf{f}}_2(\bm{\Phi}_2),
\end{align}
where $\widetilde{S}: = \widetilde{G}_1 \widetilde{M}_1^{-1} \widetilde{G}_1^T + G_2 M_2^{-1} G_2^T$ is the Schur complement of the upper $2\times 2$ block of the ROM-FEM monolithic problem \eqref{AdC:ROMsystem}. At this juncture, we point out that we may safely expect the matrix $\widetilde{M}_1: = \widetilde{U}^T_0 M_1 \widetilde{U}_0$ to be invertible because $M_1$ is a symmetric, positive definite matrix, and multiplication by the orthogonal matrix $\widetilde{U}_0$ preserves the rank of the matrix. Now, the system (\ref{AdC:ROMsystem}) can be equivalently written as: 
\begin{align}
\begin{split}
\dot{y} &= f(t,y,z) \\
0 &= g(t,y,z)
\end{split}
\end{align}

Extension of the IVR scheme to the ROM-FEM system \eqref{AdC:ROMsystem} requires the Jacobian $\partial_z g = \widetilde{S}$ to be non-singular for all $t$. In the case of the FEM-FEM coupled system \eqref{AdC:coupled425} conditions on the Lagrange multiplier space were given in \citeauthorandyear{AdC:CAMWA} that correspond to properties of the matrices $G_1, G_2$, and ensure that the FEM-FEM Schur complement is symmetric and positive definite. In the case of the ROM-FEM coupled problem we have observed numerically that the corresponding Schur complement $\widetilde{S}$ is nonsingular. A formal proof and a sufficient condition for $\widetilde{S}$ to be symmetric and positive definite is in progress and will be reported in a forthcoming paper. 

The ROM-FEM monolithic system \eqref{AdC:ROMsystem} is the basis for the new hybrid partitioned IVR scheme. Although structurally, this problem is similar to the monolithic system  \eqref{AdC:matrixfull} underpinning the FEM-FEM scheme, there are some algorithmic distinctions that 
we wish to highlight. Most notably, the  partitioned ROM-FEM IVR algorithm has two phases: an offline phase to compute the ROM and an online phase where the ROM is used in the partitioned scheme to solve the coupled system.
For example, in the context of a PDE-constrained optimization algorithm that requires multiple solutions of the coupled problem, the first phase would be conducted offline before the optimization loop, and then the second phase would run at each optimization iteration.

\medskip
\noindent
\textbf{Computation of the reduced order model (Offline)}
\begin{enumerate}
\item 
	Use an appropriate FOM to simulate the solution on $\Omega_1$ and collect samples for the snapshot matrix $X$.
	Compute the SVD of the adjusted snapshot matrix $X_0$ containing zeros on all Dirichlet rows of $X$:
	$X_0 = U_0 \Sigma_0 V^T_0$.
\item Given a threshold $\delta >0$, define the reduced basis matrix $\widetilde{U}_0$ by discarding all columns in $U_0$ corresponding to singular values less than $\delta$.
\item Precompute the ROM matrices:
$$
\widetilde{M_1}:=\widetilde{U}_0^T M_1\widetilde{U}_0;\quad 
\widetilde{D}_1:=\widetilde{U}_0^T D_1 \widetilde{U}_0;\quad
\widetilde{A}_1:= \widetilde{U}_0^T A_1 \widetilde{U}_0;\quad\mbox{and}\quad \widetilde{G}_1:= G_1 \widetilde{U}_0\,.
$$
\end{enumerate} 

\medskip
\noindent
\textbf{Solution of the coupled ROM-FEM system for $t \in [0,T]$ (Online)}
\begin{enumerate}
\vspace{1mm}
\item Choose an explicit time integration scheme, i.e., the operator $D^n_t(\boldsymbol{\varphi})$.
\item For $n=0,1,\ldots$ use $\bm{\varphi}_R^n$ to compute the vector 
$$
\widetilde{\bm{f}}_1^n :=\widetilde{U}_0^T \mathbf{f}_1 - (\widetilde{D}_1 + \widetilde{A}_1)\bm{\varphi}_R^n - \widetilde{U}^T_0(D_1+A_1)\bm{\beta}.
$$
\item Use $\bm{\Phi}_2^n$ to compute the vector $\widetilde{\boldsymbol{f}}_2^n := \overline{\mathbf{f}}_2(\bm{\Phi}_2^n) = \mathbf{f}_2 - (D_2 + A_2) \bm{\Phi}_2^n
$
\item Solve the Schur complement system 
$$
\big( \widetilde{G}_1 \widetilde{M}_1^{-1} \widetilde{G}_1^T  + 
G_2 M_2^{-1} G_2^T \big) \bm{\lambda}^n = 
\widetilde{G}_1 \widetilde{M}_1^{-1} \widetilde{U}_0^T\widetilde{\bm{f}}_1^n -  G_2 M_2^{-1} \widetilde{\boldsymbol{f}}_2^n
$$
for $\bm{\lambda}^n$. Compute $\widetilde{G}_1^T \bm{\lambda}^n$ and $G_2^T \bm{\lambda}^n$.
 \item Solve the system $\widetilde{M}_1 D^n_t(\bm{\varphi}_R) = \widetilde{\bm{f}}_1^n - \widetilde{G}_1^T \bm{\lambda}^n$ 
	 and project the ROM solution to the state space of the full order model:
 $
 \bm{\Phi_1} :=\widetilde{U}_{0}\bm{\varphi}_R + \bm{\beta}_1; \quad
$
 \vspace{.25mm}
 \item Solve the system $M_2 D^n_t(\bm{\Phi}_2) = \widetilde{\boldsymbol{f}}_2^n + G_2^T \bm{\lambda}^n$.
\end{enumerate}

\subsection{IVR EXTENSION TO ROM-ROM COUPLING} \label{sec:ROM-ROM}
In this section we briefly explain the extension of the IVR scheme to a ROM-ROM case, i.e., when a ROM on $\Omega_1$ is coupled to another  ROM on $\Omega_2$.
For $j=1,2$, let $\widetilde{U}_{j,0}$ and $\bm{\beta}_{k,j}$ be the reduced basis matrix and the vectors \eqref{eq:beta} constructed on $\Omega_j$ according to the workflow in Section \ref{ROMsec}. We note here that our framework does not require the two ROMs being coupled
to have the same number of reduced basis modes. 
For simplicity, we shall assume again time-independent Dirichlet boundary conditions, so that the vectors $\bm{\beta}_{k,j}$ reduce to a vector $\bm{\beta}_j$ whose Dirichlet coefficients equal the nodal values of $g(x)$ and the free coefficients are zero.

As in Section \ref{AdC:sec:ROMsystem}, we implement the second stage of the POD-based model order reduction directly in the transformed semi-discrete monolithic problem \eqref{AdC:matrixfull}. Specifically, we perform a state transformation of both subdomain equations using the ansatz $\bm{\Phi}_{1} = \widetilde{U}_{1,0} \bm{\varphi}_{R} + \bm{\beta}_1$ for the first equation, and the ansatz
$\bm{\Phi}_2 = \widetilde{U}_{2,0} \bm{\psi}_R + \bm{\beta}_2$ for the second equation.
Then, we multiply the first equation by $\widetilde{U}_{1,0}^T$ and the second equation by $\widetilde{U}_{2,0}^T$.  The resulting ROM-ROM monolithic system is the basis for the ROM-ROM partitioned IVR algorithm which we state below.

\medskip
\noindent
\textbf{Computation of the reduced order models (Offline)}
\begin{enumerate}
\item 
	For $j=1,2$, use an appropriate FOM to simulate the solution on $\Omega_j$ and collect samples for the snapshot matrix $X_j$.
        Compute the SVD of the adjusted snapshot matrix $X_{j,0} = U_{j,0} \Sigma_{j,0} V^T_{j,0}$.
\item Given a threshold $\delta_j >0$, define the reduced basis matrices $\widetilde{U}_{j,0}$ by discarding all columns in $U_{j,0}$ corresponding to singular values less than $\delta_j$ for $j=1,2$.
\item For $j=1,2$, precompute the ROM matrices:
$$
\widetilde{M_j}:=\widetilde{U}_{j,0}^T M_j\widetilde{U}_{j,0};\quad 
\widetilde{D}_j:=\widetilde{U}_{j,0}^T D_j \widetilde{U}_{j,0};\quad
\widetilde{A}_j:= \widetilde{U}_{j,0}^T A_j \widetilde{U}_{j,0};\quad\mbox{and}\quad \widetilde{G}_j:= G_j \widetilde{U}_{j,0}\,.
$$
\end{enumerate} 

\medskip
\noindent
\textbf{Solution of the coupled ROM-ROM system for $t \in [0,T]$ (Online)}
\begin{enumerate}
\vspace{1mm}
\item Choose an explicit time integration scheme for each subdomain, i.e., an operator $D^n_{j,t}(\boldsymbol{\varphi})$, $j=1,2$.
\item For $n=0,1,\ldots$ use $\bm{\varphi}_R^n$ to compute the vector 
$$
\widetilde{\bm{f}}_1^n :=\widetilde{U}_{1,0}^T \mathbf{f}_1 - (\widetilde{D}_1 + \widetilde{A}_1)\bm{\varphi}_R^n - \widetilde{U}^T_{1,0}(D_1+A_1)\bm{\beta}_1.
$$
\item For $n=0,1,\ldots$ use $\bm{\psi}_R^n$ to compute the vector 
$$
\widetilde{\bm{f}}_2^n :=\widetilde{U}_{2,0}^T \mathbf{f}_2 - (\widetilde{D}_2 + \widetilde{A}_2)\bm{\psi}_R^n - \widetilde{U}^T_{2,0}(D_2+A_2)\bm{\beta}_2.
$$
\item Solve the Schur complement system 
$$
\big( \widetilde{G}_1 \widetilde{M}_1^{-1} \widetilde{G}_1^T  + 
\widetilde{G}_2 \widetilde{M}_2^{-1} \widetilde{G}_2^T \big) \bm{\lambda}^n = 
\widetilde{G}_1 \widetilde{M}_1^{-1} \widetilde{U}_{1,0}^T\widetilde{\bm{f}}_1^n -  
\widetilde{G}_2 \widetilde{M}_2^{-1}  \widetilde{U}_{2,0}^T\widetilde{\boldsymbol{f}}_2^n
$$
for $\bm{\lambda}^n$. Compute $\widetilde{G}_1^T \bm{\lambda}^n$ and $\widetilde{G}_2^T \bm{\lambda}^n$.
 \item Solve the system $\widetilde{M}_1 D^n_{1,t}(\bm{\varphi}_R) = \widetilde{\bm{f}}_1^n - \widetilde{G}_1^T \bm{\lambda}^n$.
 \vspace{.25mm}
 \item Solve the system $\widetilde{M}_2 D^n_{2,t}(\bm{\psi}_R) = \widetilde{\bm{f}}_2^n + \widetilde{G}_2^T \bm{\lambda}^n$.
 \item Project the ROM solutions $\bm{\varphi}_R, \bm{\psi}_R$ to the state spaces of the full order models on $\Omega_1$ and $\Omega_2$:
 $$
 \bm{\Phi}_1 := \widetilde{U}_{1,0} \bm{\varphi}_R + \bm{\beta}_1; \quad
 \bm{\Phi}_2 := \widetilde{U}_{2,0} \bm{\psi}_R + \bm{\beta}_2.
 $$
\end{enumerate}

\section{NUMERICAL EXAMPLES} \label{sec:numerical}

\begin{figure}[!ht]
\begin{center}
	\subfigure[Initial conditions]{\includegraphics[scale=.13]{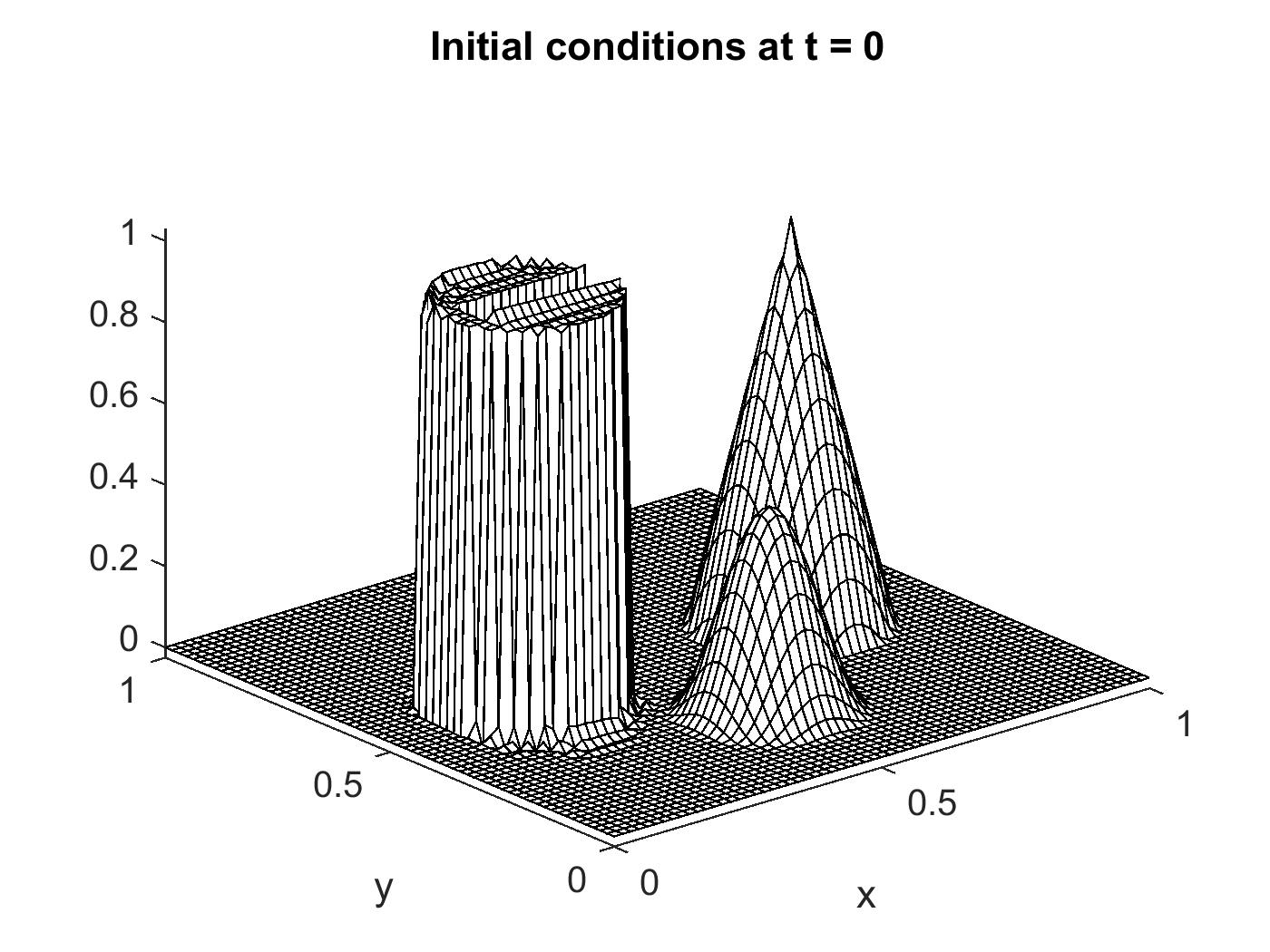}\label{AdC:fig:ICS}} \subfigure[Meshes used to discretize $\Omega_1$ (blue) and $\Omega_2$ (red)]{\includegraphics[scale=.12]{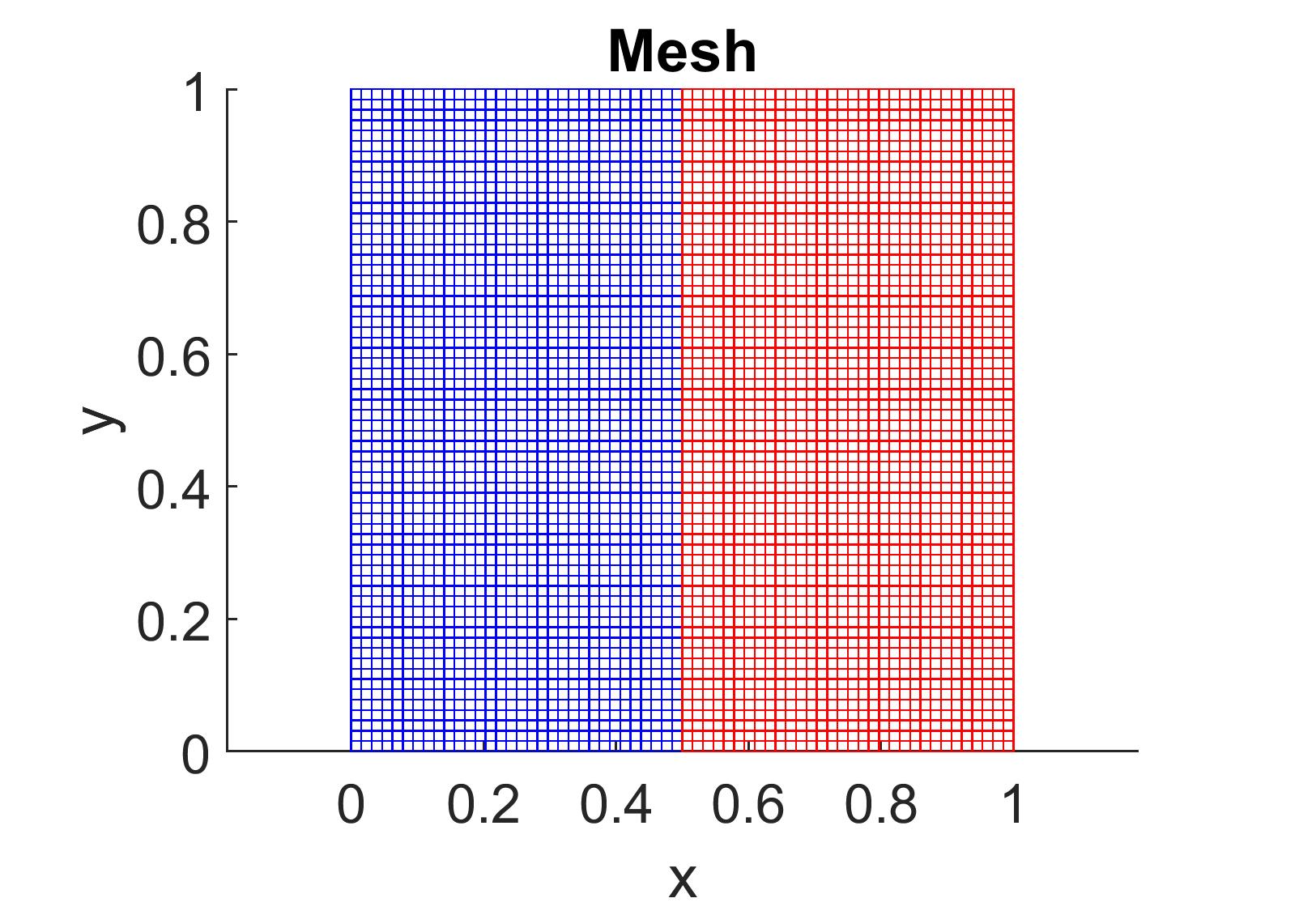}\label{AdC:fig:unifMesh}} \caption{Initial conditions and domain decomposition/mesh 
	for our model 2D transmission problem}
\end{center}
\end{figure}

To evaluate our schemes, we adapt the solid body rotation test for \eqref{AdC:strongForm}
from \citeauthorandyear{Leveque_96_SINUM}.  The problem is posed on the unit square $\Omega = (0,1)\times (0,1)$ and
the following rotating advection field $(0.5 - y, x - 0.5)$ is specified. The initial conditions for this test problem comprise a cone, cylinder, and a smooth hump, and are shown in Figure \ref{AdC:fig:ICS}.
We impose homogeneous Dirichlet boundary conditions on the non-interface boundaries $\Gamma_i$, $i=1,2$.  
We consider herein two problem configurations for \eqref{AdC:strongForm} that differ in the choice of the diffusion coefficient. The ``pure advection" case corresponds to $\kappa_i = 0$, 
and the ``high Pecl\'{e}t" case corresponds to $\kappa_i = 10^{-5}$.  
In the former case we adjust the boundary condition so that the boundary values are specified only on the inflow parts of $\Gamma_i$.
In all our tests, we run the simulations for one full rotation, i.e., the final simulation time is set to $t = 2\pi$.
It can be shown that, for the pure advection variant of this problem, the solution at the final time $t=2\pi$
should be the same as the initial solution \citeauthorandyear{Leveque_96_SINUM}.

%

Suppose $\Omega$ is divided in half vertically by the line $x = 0.5$, and let $\Omega_1$ 
and $\Omega_2$ denote the left and right side of the domain, respectively, as shown in Figure \ref{AdC:fig:unifMesh}. Let $\gamma$ denote the interface $(x = 0.5)$ between the two sides, and let $\Gamma_i = \partial\Omega_i \backslash \gamma$ for $i = 1,2$. We take $\mathbf{n}_\gamma$ to be the unit normal on the interface pointing toward $\Omega_2$. 
In this section, we present select results for solving the model advection-diffusion interface problem \eqref{AdC:strongForm} by 
performing both ROM-FEM and ROM-ROM coupling in the two subdomains, $\Omega_1$ and $\Omega_2$.  
The coupled ROM-FEM and ROM-ROM problems are solved by using the IVR partitioned schemes formulated in Section \ref{AdC:sec:ROMsystem} and Section \ref{sec:ROM-ROM}, respectively.
We compare our ROM-FEM and ROM-ROM solutions to results obtained by employing our IVR partitioned scheme to perform FEM-FEM coupling 
between the two subdomains (see Section \ref{AdC:FEM scheme}).  For comparison purposes, we also include results obtained by building
a global (uncoupled) FEM model as well as a global ROM in the full domain $\Omega$.
%

For the FEM discretizations, we employ a uniform spatial resolution of $\frac{1}{64}$ in both the $x$ and $y$ directions.
 The ROMs are developed from snapshots collected from a monolithic FEM 
 discretization of $\Omega$ using the approach described in Section \ref{ROMsec},
 with snapshots collected at intervals $\Delta_s t = 1.35 \times 10^{-2}$ 
 and $\Delta_s t = 6.73 \times 10^{-3}$ for the pure advection and high Pecl\'{e}t 
 variants of our test case, respectively.  These snapshot selection
 strategies yield 466 and 933 snapshots for the two problem variants, 
 respectively.  
 All ROMs evaluated herein are run in the reproductive regime, that is, 
 with the same parameter values, boundary conditions and initial conditions
 as those used to generate the snapshot set from which these models were constructed;
 predictive ROM simulations will be considered in a subsequent publication.
In general, between 20-25 modes 
are needed to capture 90\% of the snapshot energy and between 50-65 modes are needed 
to capture 99.999\% of the snapshot energy for both problem variants, 
where the snapshot energy fraction is defined as $1-\delta$. 
As noted in Section \ref{sec:ROM-ROM}, for the ROM-ROM couplings, we allow the bases in $\Omega_1$ 
and $\Omega_2$ to have different numbers of modes,
denoted by $N_{R, \text{left}}$ and $N_{R, \text{right}}$, respectively.  Hence, the number of modes
required to capture a given snapshot energy fraction varies slightly between
the two subdomains.
 All simulations are performed using an explicit 
 4th order Runge-Kutta (RK4) scheme with time-step $\Delta t = 3.37\times 10^{-3}$,
 the time-step computed by the Courant-Friedrichs-Lewy (CFL) condition for this problem.

 In the results below, we report for the various models evaluated the following relative errors as a function of the 
 basis size and the total online CPU time:
 \begin{equation} \label{eq:err1}
	 \epsilon := \frac{||X_{2\pi} - F_{2\pi}||_2}{||F_{2\pi}||_2}.
 \end{equation}
 In \eqref{eq:err1}, $X \in \{R, FF, RF, RR\}$, where
$R$ 
 denotes the global ROM solution computed in all of $\Omega$, $FF$ denotes a FEM-FEM 
 coupled solution, $RF$ denotes a ROM-FEM coupled solution, and $RR$ denotes a 
 ROM-ROM coupled solution.  
 The subscripts in \eqref{eq:err1} denote the time at 
 which a given solution is evaluated, i.e., $RF_{2\pi}$ is the ROM-FEM solution
 at time $t=2\pi$.  The reference solution in \eqref{eq:err1}, denoted 
 by $F_{2\pi}$, is the 
global FEM solution computed in all of $\Omega$ at time $t=2\pi$.
 For the pure advection problem, we additionally report:
 \begin{equation} \label{eq:err2}
	 \epsilon_0 := \frac{||X_{0} - X_{2\pi}||_2}{||X_{2\pi}||_2},
 \end{equation}
 for $X \in \{F, R, FF, RF, RR\}$.  As shown in \citeauthorandyear{Leveque_96_SINUM}, for the exact solution to the 
 pure advection problem, $\epsilon_0$ is identically zero.  



\begin{figure}[!ht]
	\begin{center}
		\subfigure[Pure Advection]{\includegraphics[scale=0.33]{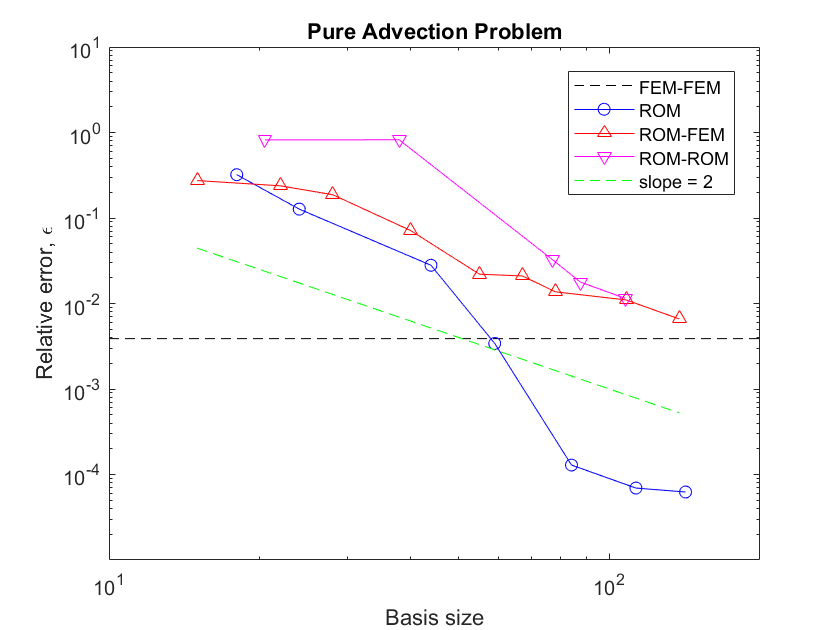}} \subfigure[High Pecl\'{e}t]{\includegraphics[scale=0.33]{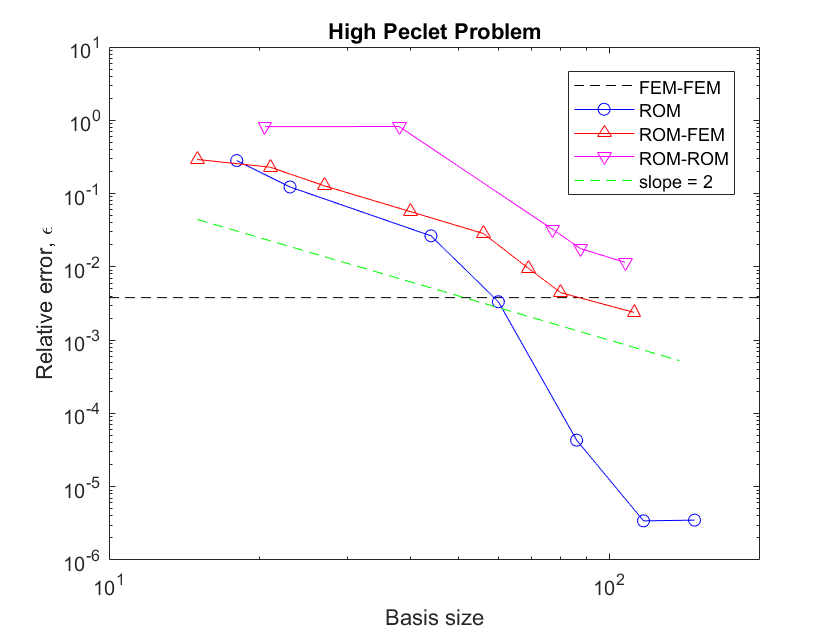}}\caption{Relative errors \eqref{eq:err1} with respect to the global FEM solution 
		as a function of the POD basis size for different discretizations of the pure advection (a) and high Pecl\'{e}t (b) variants of 
		our model transmission problem.} \label{fig:conv}
	\end{center}
	\end{figure}

First, in Figure \ref{fig:conv}, we plot the relative error $\epsilon$ in \eqref{eq:err1} as a function
of the POD basis size for the various couplings and the two problem variants considered herein. All errors are calculated
with respect to the global FEM solution computed in all of $\Omega$.  
For the ROM-ROM couplings, the basis size in Figure \ref{fig:conv} is obtained by calculating the average of the
basis sizes in $\Omega_1$ and $\Omega_2$, denoted by $N_{R, \text{left}}$ and $N_{R, \text{right}}$ respectively.  The reader can observe that all models 
exhibit convergence with respect to the basis size.  In particular, the ROM-FEM and ROM-ROM solutions 
converge at a rate of approximately two.  For the pure advection problems, the ROM-FEM and ROM-ROM
solutions appear to be approaching the FEM-FEM error with basis refinement, and the ROM-ROM 
solution appears to be converging to the ROM-FEM solution.  It is interesting 
to observe that the global ROM solutions achieve a greater accuracy than the FEM-FEM coupled solutions.
Moreover, for the high Pecl\'{e}t version of the problem, the ROM-FEM coupled solution can achieve 
an accuracy that is slightly better than the FEM-FEM coupled solution.
This behavior is likely due to the fact that the ROM solution was created using snapshots from 
a global FEM solution, which is more accurate than the coupled FEM-FEM solution. 

\begin{figure}
	\begin{center}
		\subfigure[Pure Advection]{\includegraphics[scale=0.18]{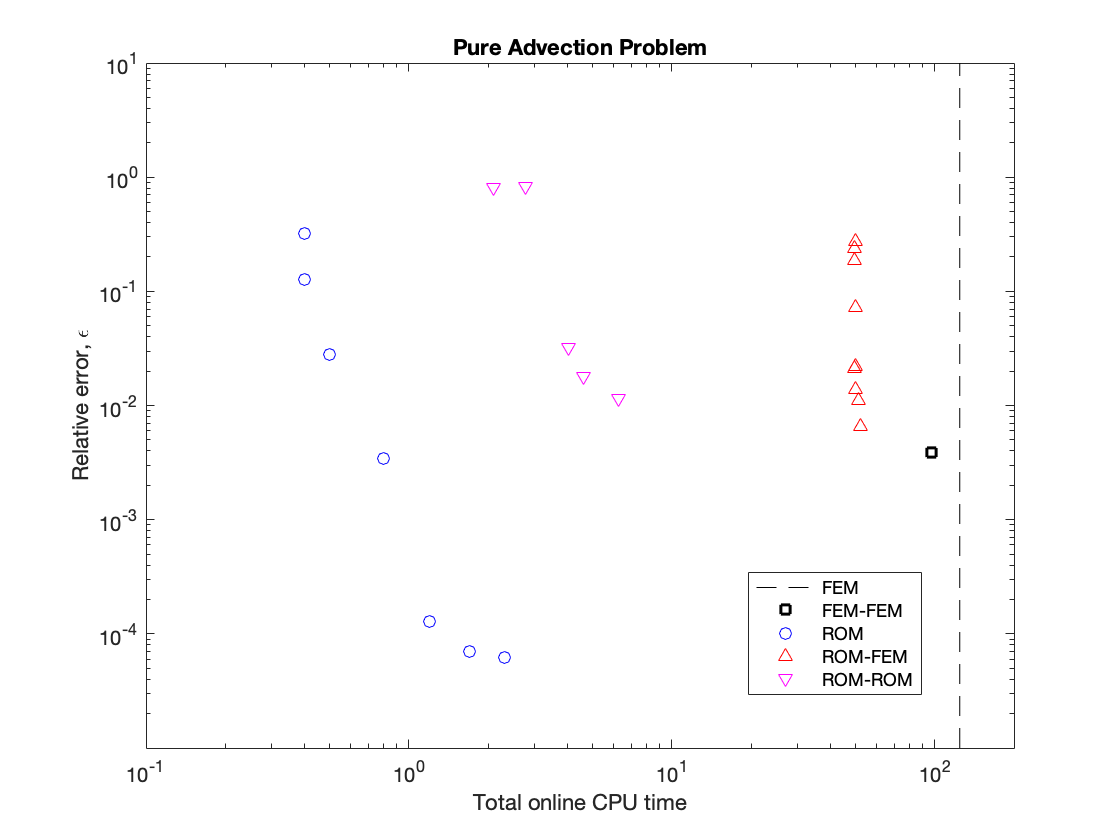}} \subfigure[High Pecl\'{e}t]{\includegraphics[scale=0.18]{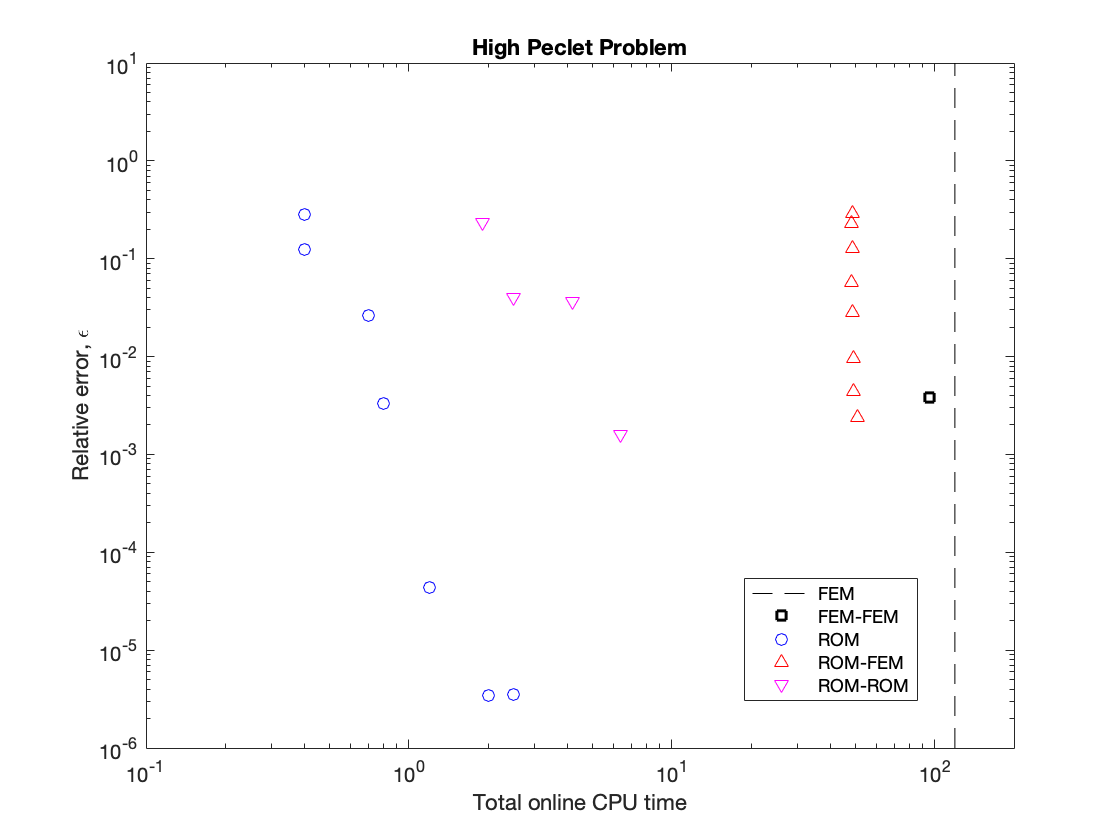}}\caption{Pareto plot (relative errors \eqref{eq:err1} as a function of the total online CPU time) for different discretizations of the pure advection (a) and high Pecl\'{e}t (b) variants of our model transmission problem.} \label{fig:pareto}
	\end{center}
	\end{figure}

In evaluating the viability of a reduced model, it is important to consider not only the model's accuracy,
but also its efficiency. Toward this effect, Figures \ref{fig:pareto}(a) and (b) show Pareto plots for the models evaluated 
on the pure advection and high Pecl\'{e}t problems, respectively.  In these figures, we plot 
the relative errors \eqref{eq:err1} as a function of the total online CPU time.  As expected, the global FEM
and FEM-FEM models require the largest CPU time, followed by the ROM-FEM models, the ROM-ROM models and 
the global ROM models.  It is interesting to remark that the FEM-FEM discretizations 
are actually slightly faster than the global FEM discretizations.  This suggests that, in the 
case of high-fidelity models, our proposed coupling approach does not introduce any significant overhead.
While the global ROM 
achieves the most accurate solution in the shortest amount of time, we are targeting here the scenario 
where the analyst does not have access to a single domain solver, and is forced to couple
models calculated independently in different parts of the computational domain.  
The results in Figure \ref{fig:pareto} show that, by introducing ROM-FEM and ROM-ROM coupling, one can reduce 
the CPU time by 1-1.5 orders of magnitude without sacrificing accuracy.  

\begin{wrapfigure}{R!HT}{0.47\textwidth}
  \begin{center}
    \includegraphics[width=0.45\textwidth]{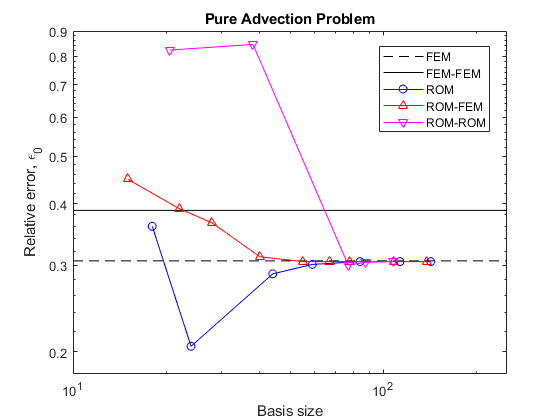}
  \end{center}
		\caption{Relative errors \eqref{eq:err2} as a function of the POD basis size for 
		different discretizations of the pure advection problem.} \label{fig:errs0}
\end{wrapfigure}

Turning our attention now to the pure advection problem, 
we plot in Figure \ref{fig:errs0} the relative errors $\epsilon_0$ in \eqref{eq:err2} as a function
of the basis size.  Again, the global FEM model is the most 
accurate, followed by the FEM-FEM, the ROM-FEM and the ROM-ROM models.  
It is interesting to observe that the global ROM surpasses the global FEM 
solution when it comes to accuracy for certain (intermediate) basis sizes.  
The primary takeaway 
from Figure \ref{fig:errs0} is that the ROM-FEM, the ROM-ROM and the global ROM solutions
asymptotically approach the global FEM solution as the basis size is refined.
This provides further verification for the models evaluated, in particular, for our
new IVR coupling approach.

\begin{figure}[!ht]
	\begin{center}
		\subfigure[Global FEM]{\includegraphics[scale = 0.12]{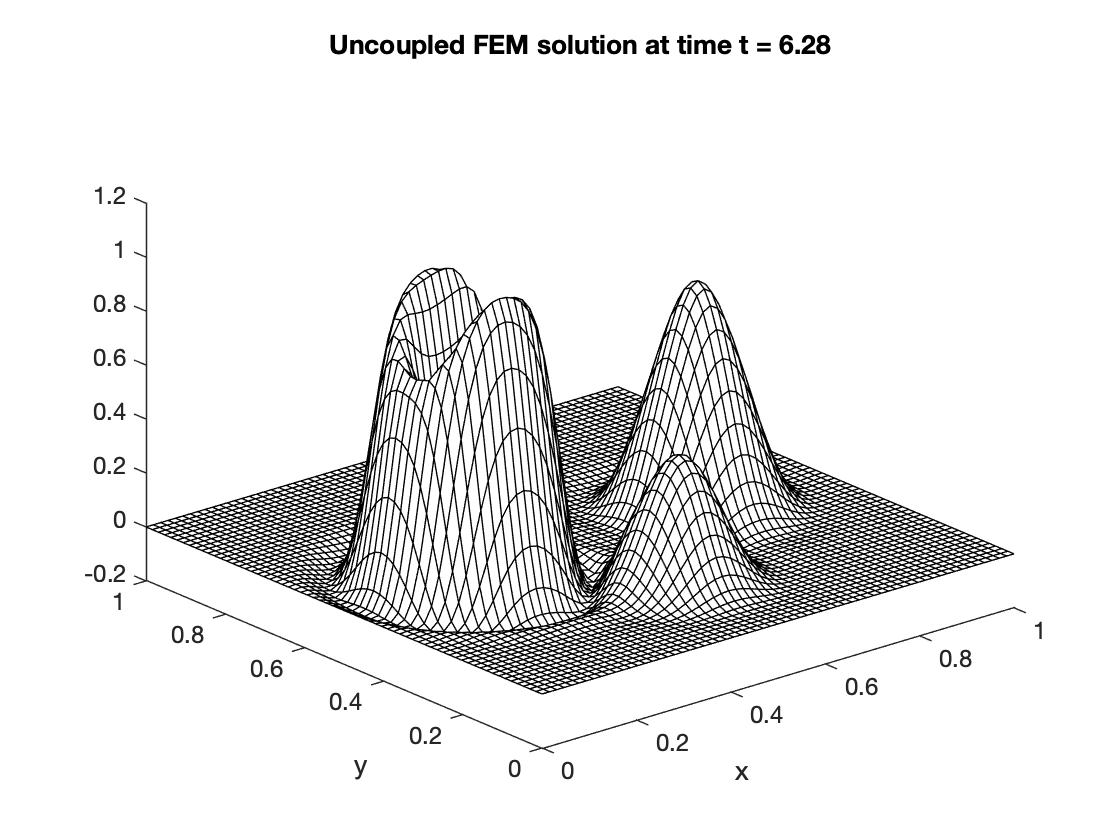}} 
		\subfigure[ROM-FEM ($N_R=80$)]{\includegraphics[scale=0.12]{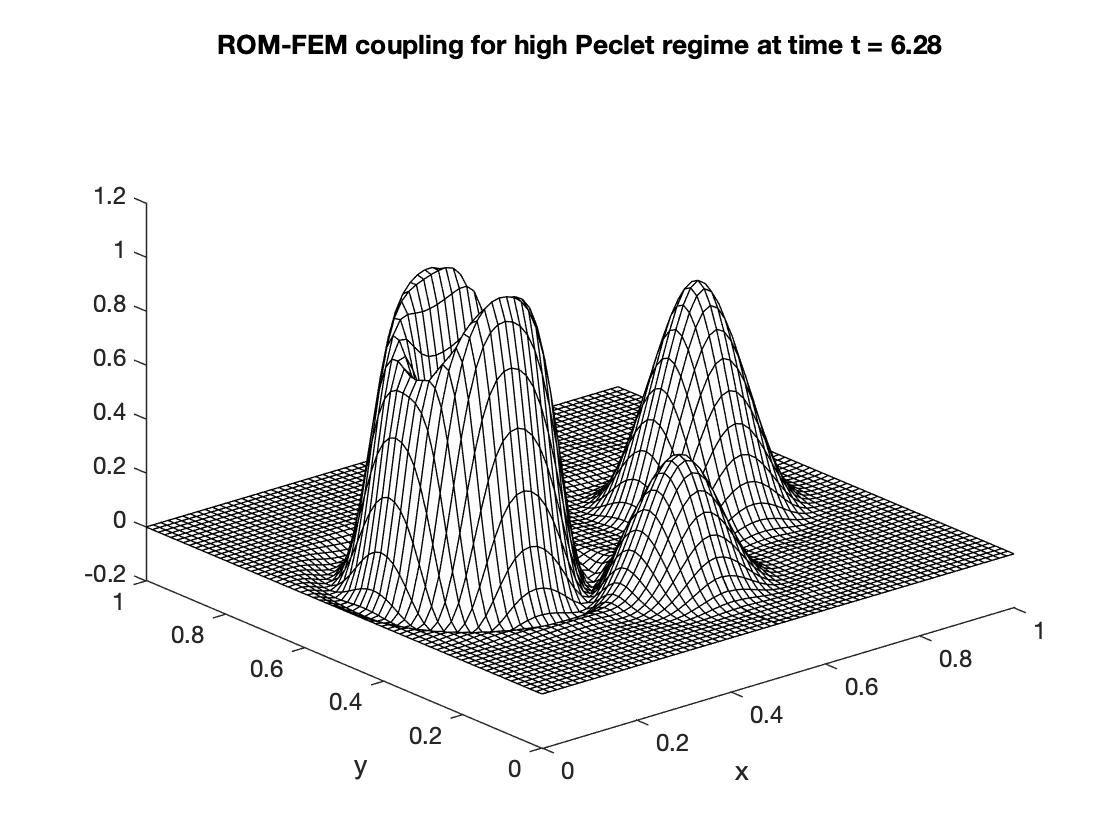}}
		\subfigure[ROM-ROM ($N_{R, \text{left}} = 112$, $N_{R, \text{right}} = 110$)]{\includegraphics[scale=0.12]{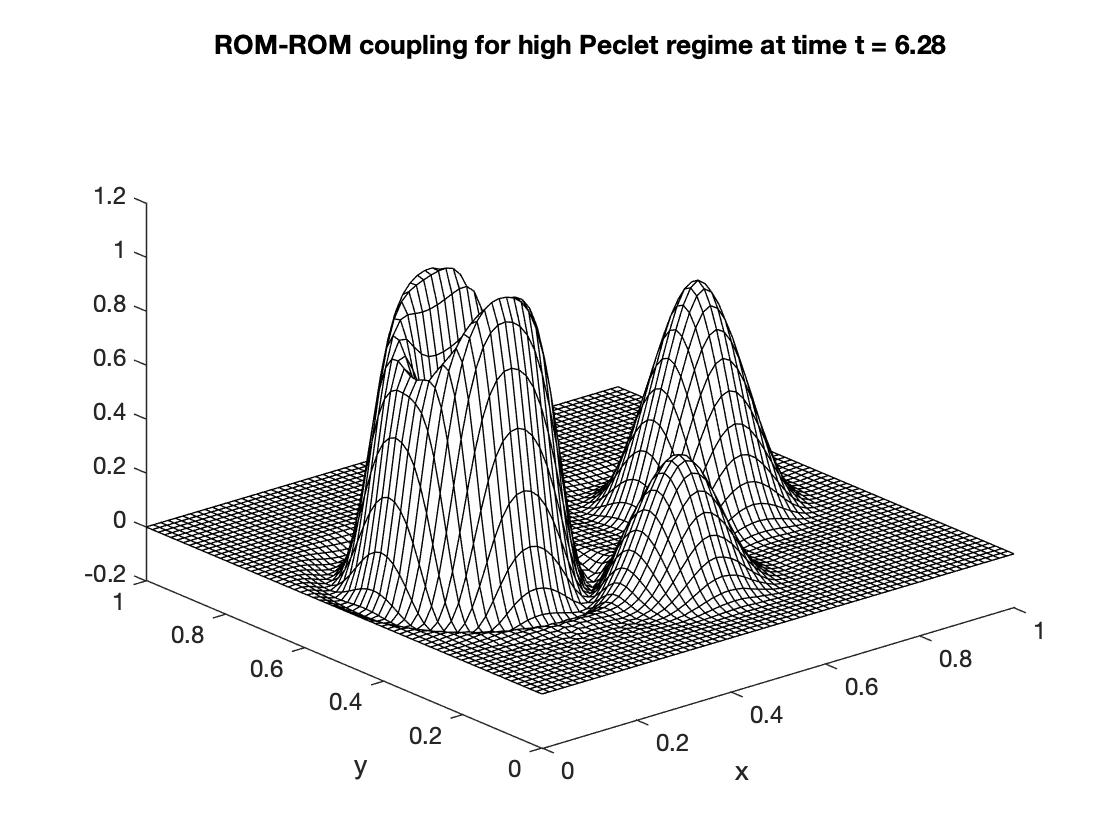}}
		\caption{Comparison of global FEM, ROM-FEM and  ROM-ROM solutions for the high Pecl\'{e}t variant 
		of our model transmission problem at the final simulation time $t=2\pi$.} \label{AdC:fig:highPe}
	\end{center}
	\end{figure}

\begin{figure}[!ht]
	\begin{center}
		\subfigure[ROM-FEM ($N_R=80$)]{\includegraphics[scale=0.19]{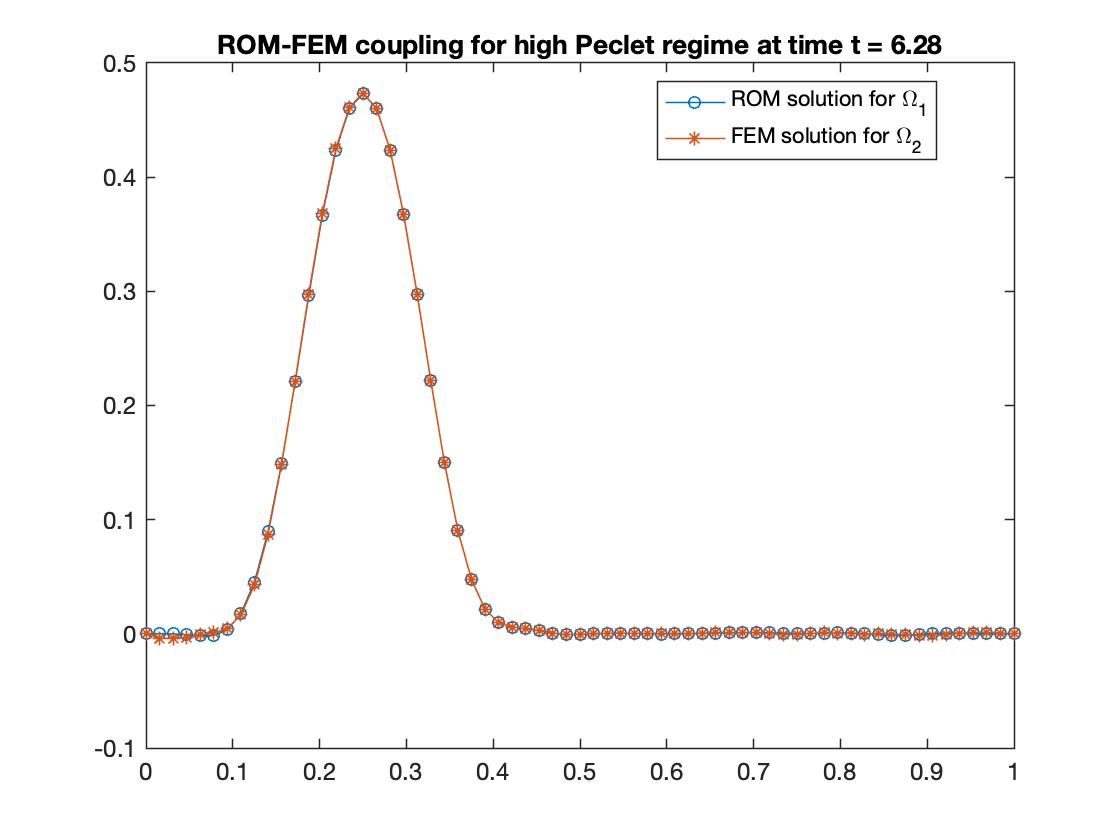}}\subfigure[ROM-ROM ($N_{R, \text{left}} = 112$, $N_{R, \text{right}} = 110$)]{\includegraphics[scale=0.19]{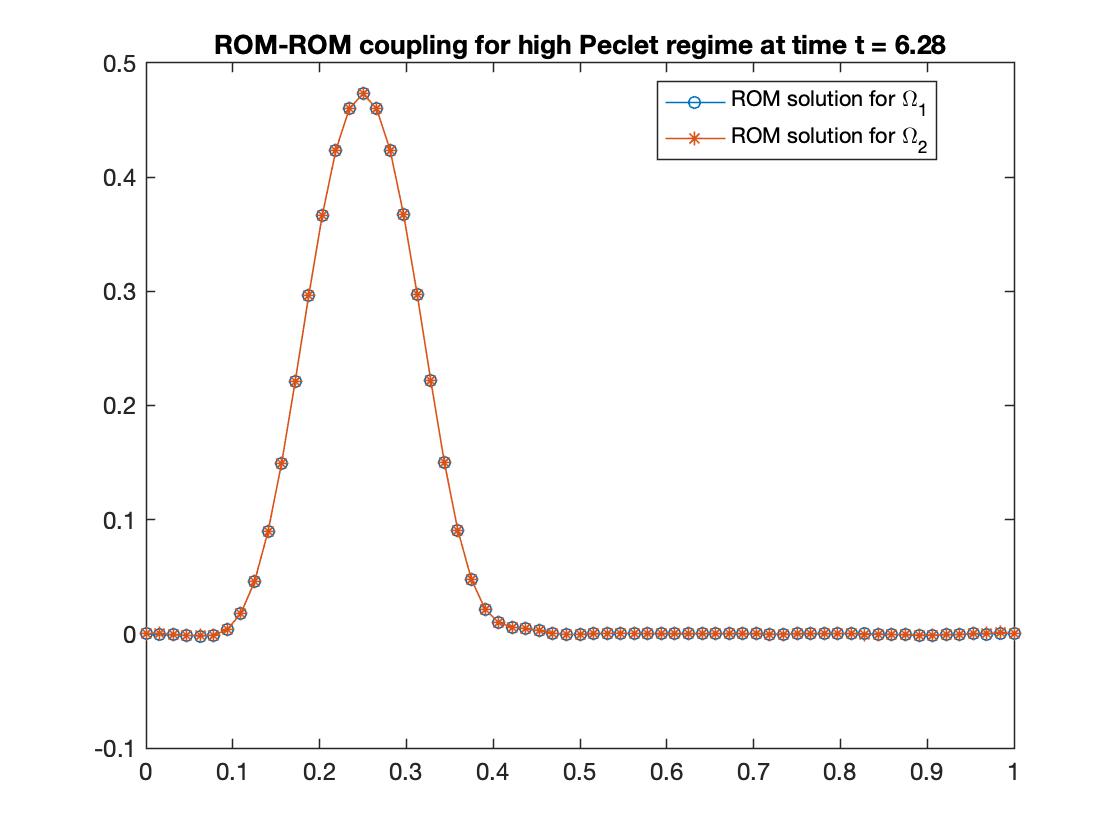}}
		\caption{Comparison of the interface  ROM-FEM and ROM-ROM solutions for high Pecl\'{e}t variant 
		of our model transmission problem at the final simulation time $t=2\pi$.} \label{AdC:fig:highPe_inter}
	\end{center}
	\end{figure}

Next, in Figure \ref{AdC:fig:highPe}, we plot some representative ROM-FEM and ROM-ROM solutions to the high Pecl\'{e}t variant
of the targeted problem at the final simulation time $2\pi$.  Also plotted is the single domain global FEM solution 
computed for this problem.  The reader can observe that all three solutions are indistinguishable 
from one another.  Figure \ref{AdC:fig:highPe_inter} plots the ROM-FEM and ROM-ROM solutions to the high Pecl\'{e}t problem
along the interface $\Gamma$ for each of the subdomains at the final simulation time $2\pi$.  
It can be seen from this figure that the solutions in $\Omega_1$ and $\Omega_2$ match incredibly well along the interface 
boundary.  This suggests that our coupling method has not introduced any spurious artifacts into the 
discretization.  We omit plots analogous to Figures \ref{AdC:fig:highPe} and \ref{AdC:fig:highPe_inter} 
for the pure advection problem for the 
sake of brevity, as they lead to similar conclusions as high Pecl\'{e}t problem results.


	\section{CONCLUSIONS} \label{sec:conc}
We presented an explicit partitioned scheme for a transmission problem that extends the approach developed in \citeauthorandyear{AdC:CAMWA} to the case of coupling a projection-based ROM 
	with a traditional finite element scheme and/or with another projection-based ROM. 
	In particular, the scheme begins with a monolithic formulation of the transmission problem and then employs a Schur complement to solve for a Lagrange multiplier representing the interface flux as a Neumann boundary condition. We constructed a ROM from a full finite element solution and then presented an algorithm to couple this reduced model with either a traditional finite element scheme or another reduced model. Our numerical results show that the ROM-FEM
	and ROM-ROM coupling produces solutions which strongly agree with those produced by a global FEM solver. 
	Additionally, implementing the ROM in one or more subdomains reduces the time and computational cost of solving the coupled system. 
	In principle, this coupling method should extend to other discretizations such as finite volume, 
	and the case of multiple ($>2$) subdomains; these scenarios will be studied in future work. Additionally, extensions to nonlinear and 
	multiphysics problems, as well as predictive runs will be considered.

	\section*{ACKNOWLEDGEMENTS}
	This work was funded by the Laboratory Directed
Research \& Development (LDRD) program at Sandia National
Laboratories, and the U.S. Department of Energy, Office of Science, Office of Advanced Scientific Computing Research under Award Number DE-SC-0000230927 and under the Collaboratory on Mathematics and Physics-Informed Learning Machines for Multiscale and Multiphysics Problems (PhILMs) project.   
The development of the ideas presented herein was funded in part by the
third author's Presidential Early Career Award for Scientists and
Engineers (PECASE).  Sandia National Laboratories is a multi-mission
laboratory managed and operated by National Technology and
Engineering Solutions of Sandia, LLC., a wholly owned subsidiary of
Honeywell International, Inc., for the U.S. Department of Energy's
National Nuclear Security Administration under contract
DE-NA0003525.  This paper describes objective technical results and
analysis. Any subjective views or opinions that might be expressed
in the paper do not necessarily represent the views of the U.S.
Department of Energy or the U.S. Government. SAND2022-7795J

\printbibliography

@article{Banks_17_JCP,
	author = {J.W. Banks and W.D. Henshaw and D.W. Schwendeman and Qi Tang},
	doi = {https://doi.org/10.1016/j.jcp.2017.01.015},
	issn = {0021-9991},
	journal = {Journal of Computational Physics},
	pages = {432 - 468},
	title = {A stable partitioned FSI algorithm for rigid bodies and incompressible flow. Part I: Model problem analysis},
	url = {http://www.sciencedirect.com/science/article/pii/S0021999117300256},
	volume = {343},
	year = {2017}}

@phdthesis{Gatzhammer_14_THESIS,
	address = {Fakultaet fuer Informatik. Informatik 5 -- Lehrstuhl fuer Wissenschaftliches Rechnen},
	author = {B. Gatzhammer},
	month = {September},
	school = {Technische Universitaet Muenchen},
	title = {Efficient and Flexible Partitioned Simulation of Fluid-Structure Interactions},
	type = {Doctoral Thesis},
	year = {2014}}

@article{Piperno_01_CMAME,
	author = {Serge Piperno and Charbel Farhat},
	doi = {http://dx.doi.org/10.1016/S0045-7825(00)00386-8},
	issn = {0045-7825},
	journal = {Computer Methods in Applied Mechanics and Engineering},
	note = {Advances in Computational Methods for Fluid-Structure Interaction},
	number = {24--25},
	pages = {3147 - 3170},
	title = {Partitioned procedures for the transient solution of coupled aeroelastic problems -- Part II: energy transfer analysis and three-dimensional applications},
	url = {http://www.sciencedirect.com/science/article/pii/S0045782500003868},
	volume = {190},
	year = {2001}}

@inproceedings{Bessette_03a_INPROC,
	author = {Greg C. Bessette and Courtenay T. Vaughan and Raymond L. Bell and Stephen W. Attaway},
	booktitle = {Proceedings of the 74 th Shock and Vibration Symposium},
	title = {Modeling Air Blast on Thin-shell Structures with {Z}apotec},
	year = {2003}}

@techreport{Randers-Pehrson_97_ARL,
	author = {Glenn Randers-Pehrson and Kenneth A. Bannister},
	institution = {Army Research Laboratory},
	month = {March},
	number = {ARL-TR-1310},
	title = {Airblast Loading Model for DYNA2D and DYNA3D},
	year = {1997}}

@article{Leveque_96_SINUM,
	author = {Randall J. LeVeque},
	doi = {10.1137/0733033},
	journal = {SIAM Journal on Numerical Analysis},
	number = {2},
	pages = {627-665},
	publisher = {SIAM},
	title = {High-Resolution Conservative Algorithms for Advection in Incompressible Flow},
	url = {http://link.aip.org/link/?SNA/33/627/1},
	volume = {33},
	year = {1996}}

@article{deBoer_07_CMAME,
	author = {A. {de Boer} and A.H. van Zuijlen and H. Bijl},
	doi = {http://dx.doi.org/10.1016/j.cma.2006.03.017},
	issn = {0045-7825},
	journal = {Computer Methods in Applied Mechanics and Engineering},
	note = {Domain Decomposition Methods: recent advances and new challenges in engineering},
	number = {8},
	pages = {1515 - 1525},
	title = {Review of coupling methods for non-matching meshes},
	url = {http://www.sciencedirect.com/science/article/pii/S0045782506002817},
	volume = {196},
	year = {2007}}

@article{Bochev_05_SIREV,
	author = {Pavel Bochev and R. B. Lehoucq},
	doi = {10.1137/S0036144503426074},
	journal = {SIAM Review},
	number = {1},
	pages = {50-66},
	publisher = {SIAM},
	title = {On the Finite Element Solution of the Pure {N}eumann Problem},
	url = {http://link.aip.org/link/?SIR/47/50/1},
	volume = {47},
	year = {2005}}

@article{AdC:CAMWA,
	author = {K. Peterson and P. Bochev and P. Kuberry},
	issue = {2},
	journal = {Computers \& Mathematics with Applications},
	pages = {459--482},
	title = {Explicit synchronous partitioned algorithms for interface problems based on Lagrange multipliers},
	volume = {78},
	year = {2019}}

@article{AdC:RINAM,
	author = {K. Chad Sockwell and Kara Peterson and Paul Kuberry and Pavel Bochev and Nat Trask},
	journal = {Results in Applied Mathematics},
	pages = {100--110},
	title = {Interface Flux Recovery coupling method for the ocean--atmosphere system},
	volume = {8},
	year = {2020}}

@phdthesis{LeGresley:2005,
	address = {Stanford, CA},
	author = {LeGresley, P.},
	school = {Stanford University},
	title = {{Application of Proper Orthogonal Decomposition (POD) to Design Decomposition Methods}},
	year = 2005}

@inproceedings{LeGresley:2003,
	author = {P. LeGresley and Alonso, J.},
	pages = {41st Aerospace Sciences Meeting and Exhibit, Reno, Nevada},
	title = {Dynamic domain decomposition and error correction for reduced order models},
	year = 2003}

@techreport{Buffoni:2007,
	author = {M. Buffoni and Telib, H. and Iollo, A.},
	institution = {INRIA Report},
	number = {6383},
	title = {{Iterative Methods for Model Reduction by Domain Decomposition}},
	year = 2007}

@inproceedings{Cinquegrana:2011,
	author = {D. Cinquegrana and A. Viviani and R. Donelli},
	pages = {AIMETA 2011- XX Congresso dell'Associazione Italiana di Meccanica Teorica e Applicata, Bologna, ITA},
	title = {{A hybrid method based on POD and domain decomposition to compute the 2-D aerodynamic flow field - incompressible validation}},
	year = 2011}

@article{Maier:2014,
	author = {I. Maier and B. Haasdonk},
	journal = {Applied Numerical Mathematics},
	pages = {31-48},
	title = {A Dirichlet-Neumann reduced basis method for homogeneous domain decomposition problems},
	volume = {78},
	year = {2014}}

@article{Ammar:2011,
	author = {A. Ammar and F. Chinestra and E. Cueto},
	journal = {International Journal for Multiscale Computational Engineering},
	pages = {17-33},
	title = {Coupling finite elements and proper generalized decompositions},
	volume = {9},
	year = {2011}}

@article{Iapichino:2016,
	journal = {Computers \& Mathematics with Applications},
	author = {L. Iapichino and A. Quarteroni and G. Rozza}, 
	number = {1},
	pages = {408-430},
	title = {Reduced basis method and domain decomposition for elliptic problems in networks and complex parametrized geometries},
	volume = {71},
	year = {2016}}

@article{Eftang:2013,
	author = {Eftang, Jens L. and Patera, Anthony T.},
	journal = {International Journal for Numerical Methods in Engineering},
	number = {5},
	pages = {269-302},
	title = {Port reduction in parametrized component static condensation: approximation and a posteriori error estimation},
	volume = {96},
	year = {2013}}

@article{Eftang:2014,
	author = {Eftang, Jens L. and Patera, Anthony T.},
	journal = {Advanced Modeling and Simulation in Engineering Sciences},
	number = {3},
	title = {{A port-reduced static condensation reduced basis element method for large component-synthesized structures: approximation and A Posteriori error estimation}},
	volume = {1},
	year = {2014}}

@article{Hoang:2021,
	author = {Chi Hoang and Youngsoo Choi and Kevin Carlberg},
	journal = {Computer Methods in Applied Mechanics and Engineering},
	pages = {113997},
	title = {Domain-decomposition least-squares Petrov--Galerkin (DD-LSPG) nonlinear model reduction},
	volume = {384},
	year = {2021}}

@article{Lucia:2003,
	author = {David J. Lucia and Paul I. King and Philip S. Beran},
	issn = {0045-7930},
	journal = {Computers \& Fluids},
	number = {7},
	pages = {917-938},
	title = {Reduced order modeling of a two-dimensional flow with moving shocks},
	volume = {32},
	year = {2003}}

@inproceedings{Lucia:2001,
	author = {D. Lucia and King, P. and Beran, P. and Oxley, M.},
	pages = {15th AIAA Computational Fluid Dynamics Conference, Anaheim, CA},
	title = {Reduced order modeling for a one-dimensional nozzle flow with moving shocks},
	year = 2001}

@article{Maday:2004,
	journal = {SIAM J. Sci. Comput.},
	author = {Y. Maday and E. Ronquist},
	number = {1},
	pages = {240-258},
	title = {The reduced basis element method: application to a thermal fin problem},
	volume = {26},
	year = {2004}}

@article{Antil:2010,
	author = {H. Antil and M. Heinkenschloss and R. Hoppe and D. Sorensen},
	journal = {Computing and Visualization in Science},
	pages = {249--264},
	title = {{Domain decomposition and model reduction for the numerical solution of PDE constrained optimization problems with localized optimization variables}},
	volume = {13},
	year = {2010}}

@article{Corigliano:2015,
	author = {Corigliano, A. and Dossi, M. and Mariani, S.},
	journal = {Computer Methods in Applied Mechanics and Engineering},
	number = {C},
	pages = {127-155},
	title = {Model Order Reduction and domain decomposition strategies for the solution of the dynamic elastic--plastic structural problem},
	volume = {290},
	year = {2015}}

@article{Corigliano:2013,
	author = {Alberto Corigliano and Martino Dossi and Stefano Mariani},
	issn = {0045-7949},
	journal = {Computers \& Structures},
	note = {Computational Fluid and Solid Mechanics 2013},
	pages = {113-127},
	title = {Domain decomposition and model order reduction methods applied to the simulation of multi-physics problems in MEMS},
	volume = {122},
	year = {2013}}

@article{Kerfriden:2013,
	author = {P. Kerfriden and O. Goury and T. Rabczuk and S.P.A. Bordas},
	issn = {0045-7825},
	journal = {Computer Methods in Applied Mechanics and Engineering},
	pages = {169-188},
	title = {A partitioned model order reduction approach to rationalise computational expenses in nonlinear fracture mechanics},
	volume = {256},
	year = {2013}}

@article{Radermacher:2014,
	author = {A. Radermacher and S. Reese},
	journal = {Computational Mechanics},
	pages = {677-687},
	title = {Model reduction in elastoplasticity: proper orthogonal decomposition combined with adaptive sub-structuring},
	volume = {54},
	year = {2014}}

@article{Baiges:2013,
	author = {Joan Baiges and Ramon Codina and Sergio R. Idelsohn},
	journal = {Computer Methods in Applied Mechanics and Engineering},
	pages = {23-42},
	title = {A domain decomposition strategy for reduced order models. Application to the incompressible Navier--Stokes equations},
	volume = {267},
	year = {2013}}

@article{Holmes:1988,
	author = {Aubry, Nadine and Holmes, Philip and Lumley, John L. and Stone, Emily},
	doi = {10.1017/S0022112088001818},
	journal = {Journal of Fluid Mechanics},
	pages = {115--173},
	publisher = {Cambridge University Press},
	title = {The dynamics of coherent structures in the wall region of a turbulent boundary layer},
	volume = {192},
	year = {1988}}

@article{Gunzburger:2007,
	author = {Max D. Gunzburger and Janet S. Peterson and John N. Shadid},
	journal = {Computer Methods in Applied Mechanics and Engineering},
	pages = {1030-1047},
	title = {Reduced-order modeling of time-dependent PDEs with multiple parameters in the boundary data},
	volume = {196},
	year = {2007}}

@article{Sirovich:1987,
	author = {L. Sirovich},
	issn = {0033569X, 15524485},
	journal = {Quarterly of Applied Mathematics},
	number = {3},
	pages = {583--590},
	publisher = {Brown University},
	title = {{Turbulence and the dynamics of coherent structures, Part III: dynamics and scaling}},
	url = {http://www.jstor.org/stable/43637459},
	volume = {45},
	year = {1987}}

@book{Holmes:1996,
	author = {P. Holmes and J. Lumpley and G. Berkooz},
	publisher = {Cambridge University Press},
	title = {Turbulence, Coherent structures, dynamical systems and symmetry},
	year = 1996}


\end{document}